\documentclass[12pt]{article}

\usepackage[T1]{fontenc}
\usepackage{lmodern}
\usepackage{amsmath,amssymb,amsthm}
\usepackage{microtype}
\usepackage{needspace}
\usepackage[colorlinks=true,allcolors=blue]{hyperref}

\makeatletter
\renewcommand{\section}{%
  \@startsection{section}{1}{\z@}%
    {-3.5ex \@plus -1ex \@minus -.2ex}%
    {2.3ex \@plus .2ex}%
    {\normalfont\fontfamily{cmr}\fontseries{bx}\Large\selectfont}}
\renewcommand{\subsection}{%
  \@startsection{subsection}{2}{\z@}%
    {-3.25ex \@plus -1ex \@minus -.2ex}%
    {1.5ex \@plus .2ex}%
    {\normalfont\fontfamily{cmr}\fontseries{bx}\large\selectfont}}
\makeatother

\theoremstyle{plain}
\newtheorem{theorem}{Theorem}[section]
\newtheorem{lemma}[theorem]{Lemma}

\theoremstyle{definition}

\theoremstyle{remark}

\def\pf{\noindent{\it Proof.} }
\allowdisplaybreaks[2]
\numberwithin{equation}{section}

\begin{document}
\begin{center}
{\large \bf An Infinite Family of Overpartition Congruences Modulo 1944}
\end{center}

\begin{center}
{Li Zhang}   \vskip 2mm

Department of Mathematics and Physics   \\
 Shijiazhuang Tiedao University, Shijiazhuang 050043, Hebei, P. R. China

\vskip 2mm

Email: zhangli\_math7@stdu.edu.cn
\end{center}

\vskip 6mm \noindent 

{\bf Abstract.} 
Let $\overline p(n)$ denote the number of overpartitions of $n$.
We prove that
$\overline p(16^\alpha(96n+76))\equiv0\pmod{1944}$
for all nonnegative integers $\alpha$ and $n$, extending the
congruence for $\overline p(96n+76)$ established by Smoot.
The key step is a three-term congruence modulo $243$ relating
$\overline p(m)$, $\overline p(16m)$, and $\overline p(256m)$
for $m\equiv4\pmod6$. We derive this relation and the two
initial congruences using theta-function dissections.
A separate theta-series argument establishes divisibility by $8$.

\noindent {\bf Keywords}: Overpartitions, generating functions, congruences, theta functions

\noindent {\bf 2020 Mathematics Subject Classification}: Primary 11P83; Secondary 05A17

\section{Introduction}
An overpartition is a partition in which the first occurrence of each distinct part may be overlined. Let $\overline{p}(n)$ denote the number of overpartitions of $n$.  By convention, we set $\overline{p}(0)=1$ and $\overline{p}(n)=0$ if $n<0$.

Corteel and Lovejoy \cite{Corteel-Lovejoy-2004} showed that the generating function of $\overline{p}(n)$ is given by
\begin{equation}\label{gf-overp}
G(q):=\sum_{n\ge0}\overline p(n)q^n
 =\frac{(-q;q)_\infty}{(q;q)_\infty}
 =\frac{f_2}{f_1^2},
\end{equation}
where
\[
(a;q)_\infty := \prod_{n=0}^{\infty} (1-aq^n), \qquad  f_k:=(q^k;q^k)_{\infty}. 
\]
Throughout, $q$ is a formal variable, and congruences between power series are interpreted coefficientwise. Since each $f_k$ has constant term $1$, it is invertible in $\mathbb Z[[q]]$. Thus every quotient of products of the $f_k$ appearing below is well defined in $\mathbb Z[[q]]$.

Xia \cite{Xia-2017} conjectured that 
\begin{equation}\label{xia-conj}
\overline{p}(96n+76)\equiv0\pmod{243} \qquad (n\ge0).
\end{equation}
Huang and Yao \cite{1-Huang-Yao-2019} proved \eqref{xia-conj}
using a theta-function parametrization. Smoot
\cite[Theorem 3.5]{Smoot-2021} subsequently
established divisibility by $1944$ on the same progression using
Radu's algorithm. 
Chen, Hou, Sun and Zhang \cite{Chen-Hou-Sun-Zhang-2016} established
several congruences relating $\overline{p}(4m)$ to $\overline{p}(m)$ modulo
powers of $2$.

In this paper, we study the behavior of this progression under repeated
multiplication by $16$. The proof has two parts. First, successive
two-dissections yield congruences for the initial progressions and a
scaling relation modulo $243$ on the residue class $4$ modulo $6$.
Second, a theta-series argument supplies the factor $8$.
The resulting family extends the case $\alpha=0$ established
by Smoot.

For compactness, put
\begin{equation}\label{eq:H}
H(q):= 
\frac{f_2^2f_3^6}{f_1^3}
 - \frac{f_3f_6^4}{f_1^4f_2^2}
 + q\frac{f_6^9}{f_2f_3^3}.
\end{equation}
Our two generating-function congruences have the same right-hand side.

\begin{theorem}\label{thm-ini-gf-96}
We have
\begin{align}
\sum_{n=0}^{\infty} \overline{p}(96n+76) q^n
&\equiv  108 H(q) \pmod{243},  \label{1-gf-96} \\
\sum_{n=0}^{\infty} \overline{p}(16(96n+76)) q^n
&\equiv  108 H(q) \pmod{243}. \label{2-gf-96}
\end{align}
\end{theorem}

\begin{theorem}\label{thm-con-rel}
For every $n\geq 0$,
\begin{align}\label{con-rel}
 \overline{p}(6n+4)+ 4\cdot\overline{p}(16(6n+4)) + 4\cdot\overline{p}(16^2(6n+4)) 
\equiv 0 \pmod{243}.
\end{align}
\end{theorem}

\begin{theorem}\label{thm-inf-cong-96-m1944}
For all integers $\alpha, n\geq0$,
\begin{align}\label{inf-cong-96-m1944}
\overline{p}\bigl(16^{\alpha}(96n+76)\bigr) \equiv 0 \pmod{1944}.
\end{align}
\end{theorem}

The expression $H(q)$ vanishes modulo $9$ by \eqref{f2f3^6-f1^2f6^3}. 
Thus Theorem~\ref{thm-ini-gf-96} gives the two initial congruences modulo $243$.
Theorem~\ref{thm-con-rel} propagates them, and the modulo $8$ argument in Section~6 completes the proof of Theorem~\ref{thm-inf-cong-96-m1944}.
Section~2 collects the product identities used in the dissections. 
Sections~3 and~4 establish the two generating-function congruences in
Theorem~\ref{thm-ini-gf-96}; 
Section~5 proves Theorem~\ref{thm-con-rel}.

\section{Preliminary identities and congruences}
We use Ramanujan's theta functions in the following product forms: 
\begin{align}
 \psi(q)  &:= \sum_{n\ge0}q^{n(n+1)/2} =\frac{f_2^2}{f_1}, \label{psi} \\
\varphi(q) &:= \sum_{n=-\infty}^{\infty}q^{n^2} = \frac{f_2^5}{f_1^2f_4^2}. \label{varphi} 
\end{align}
Replacing $q$ by $-q$  in \eqref{varphi}  gives
\begin{align}
 \varphi(-q) &:=\sum_{n=-\infty}^{\infty}(-1)^nq^{n^2}= \frac{f_1^2}{f_2}. \label{varphi-q}
\end{align}

\begin{lemma}\cite[Lemma 2.3]{Jin-Zhang-2025}  \label{f-3-r}
For positive integers $k$ and $r$,
\begin{align*}
 f_{k}^{3^{r}} &\equiv f_{3k}^{3^{r-1}} \pmod{3^r}.
\end{align*}
\end{lemma}

\begin{lemma}\cite[Lemmas 2.2--2.6]{Yao-Xia-2013}; \cite[Eq. (30.12.3)]{Hirschhorn-2017}
The following 2-dissections hold.
\begin{align}
f_1^2 &= \frac{f_2f_8^5}{f_4^2f_{16}^2} - 2q\frac{f_2f_{16}^2}{f_8},\label{f1-2} \\
f_1^4 &= \frac{f_4^{10}}{f_2^2f_{8}^4} - 4q\frac{f_2^2f_{8}^4}{f_4^2},\label{f1-4} \\
\frac{1}{f_1^2} &= \frac{f_8^5}{f_2^5f_{16}^2} + 2q\frac{f_4^2f_{16}^2}{f_2^5f_8}, \label{1-f1-2} \\
\frac{1}{f_1^4} &= \frac{f_4^{14}}{f_2^{14}f_{8}^4} + 4q\frac{f_4^2f_{8}^4}{f_2^{10}}, \label{1-f1-4} \\
f_3^2 &= \frac{f_{6}f_{24}^{5}}{f_{12}^{2}f_{48}^2} - 2q^3\frac{f_{6}f_{48}^2}{f_{24}}, \label{f3-2}\\
\frac{f_3}{f_1^3} &= \frac{f_4^6f_6^3}{f_2^9f_{12}^2}+3q\frac{f_4^2f_6f_{12}^2}{f_2^7}, \label{f3-f1-3}  \\
\frac{f_1^3}{f_3} &= \frac{f_4^3}{f_{12}}-3q\frac{f_2^2f_{12}^3}{f_4f_6^2}, \label{f1-3-f3} \\
\frac{f_3^3}{f_1} &= \frac{f_4^3f_6^2}{f_2^2f_{12}}+q\frac{f_{12}^3}{f_4}, \label{f3-3-f1}  \\
\frac{f_1}{f_3^3} &= \frac{f_2f_4^2f_{12}^2}{f_{6}^7}-q\frac{f_2^3f_{12}^6}{f_4^2f_6^9}, \label{f1-f3-3} \\
\frac{f_1}{f_3} &= \frac{f_2f_{16}f_{24}^2}{f_6^2f_8f_{48}} - q\frac{f_2f_8^2f_{12}f_{48}}{f_4f_6^2f_{16}f_{24}}, \label{f1-f3}\\
\frac{f_3}{f_1} &= \frac{f_4f_6f_{16}f_{24}^2}{f_2^2f_8f_{12}f_{48}} + q\frac{f_6f_8^2f_{48}}{f_2^2f_{16}f_{24}}, \label{f3-f1}\\
\frac{1}{f_1f_3} &= \frac{f_8^2f_{12}^5}{f_2^2f_4f_6^4f_{24}^2} + q\frac{f_4^5f_{24}^2}{f_2^4f_{6}^2f_8^2f_{12}}. \label{1-f1f3}
\end{align}
\end{lemma}

\begin{lemma}\cite[Eqs. (22), (24), (38), (34), (35)]{Zhang-2025} We have
\begin{align}
 q \frac{f_1f_6^6}{f_2^2f_3^3}
& = \frac{f_2^{6}f_{3}}{f_1^{3}f_6^{2}}
  - \frac{f_2f_3^6}{f_1^2f_6^3}, \label{q-f1-f6-6} \\
 2q\frac{f_4^8}{f_2^4}
&  \equiv \frac{f_2^{20}}{f_1^8f_4^8} - \frac{f_1^3f_2^3}{f_3f_6} \pmod{9}, \label{lem-qf4^8-f2^4-m9} \\
 \frac{f_{12}^2}{f_4^6}
&  \equiv \frac{f_{3}^2}{f_1^6}
- 6q\frac{f_2^2f_4^8}{f_6^2} \pmod{27}, \label{lem-f12-2-f4-6} \\
\frac{f_1^3f_6^2}{f_2^6f_3}
&  \equiv \frac{f_2^3f_3^2}{f_1^6f_6} - 9q\frac{f_2^{16}}{f_1^8} \pmod{27},  \label{lem-f1-3-f6-2} \\
\frac{f_1^6f_{6}^4}{f_2^{12}f_3^2}
& \equiv \frac{f_2^6f_{3}^4}{f_1^{12}f_6^2}
+ 9q\frac{f_2^{19}f_3^2}{f_1^{14}f_6} \pmod{27}.
\label{lem-f1-6-f6-4}
\end{align}
\end{lemma}

\begin{lemma} We have
\begin{align}
 \frac{f_6^3}{f_2^9}
&= 9 \frac{f_3^8}{f_1^8 f_2^5 f_6}
- 8\frac{f_3^3}{f_1^9}, \label{f6-3-f2-9} \\
2 \frac{f_4^{10}}{f_2^4 f_8^4}
&= \frac{f_2^{10}}{f_1^4 f_4^4} + \frac{f_1^4}{f_2^2}, \label{f_4^{10}-f_2^4 f_8^4} \\
8 q \frac{f_8^4}{f_4^2}
& = \frac{f_2^{10}}{f_1^4 f_4^4}  - \frac{f_1^4}{f_2^2}.   \label{8q-f_8^4-f_4^2}
\end{align}
\end{lemma}
\begin{proof}
From Berndt \cite[p. 357, Eq.(4.3)]{Ramanujan-5}, we see
\begin{align}\label{Ram-id-2}
9\frac{\varphi^3(-q^3)}{\varphi(-q)}-\frac{\varphi^3(-q)}{\varphi(-q^3)} = 8\frac{\psi^3(q)}{\psi(q^3)}.
\end{align}
Substituting \eqref{psi} and \eqref{varphi-q} into \eqref{Ram-id-2} yields 
\begin{equation}\label{psi-phi-Ram-f}
9\frac{f_2f_3^6}{f_1^2f_6^3}
- \frac{f_1^6f_6}{f_2^3f_3^2}
= 8\frac{f_2^{6}f_3}{f_1^3f_6^2}.
\end{equation}
Multiplying both sides of \eqref{psi-phi-Ram-f} by $f_3^2f_6^2/(f_1^6f_2^6)$ and rearranging yields \eqref{f6-3-f2-9}.

The identities in Berndt~\cite[p.40, Entry 25]{Ramanujan-3} give 
\begin{align}
\varphi^2(q) - \varphi^2(-q) &= 8q \psi^2(q^4), \label{Ram-III-i1} \\
\varphi^2(q) + \varphi^2(-q) &= 2 \varphi^2(q^2). \label{Ram-III-i2}
\end{align}
Applying \eqref{psi}--\eqref{varphi-q} to  \eqref{Ram-III-i2} and  \eqref{Ram-III-i1} yields \eqref{f_4^{10}-f_2^4 f_8^4} and \eqref{8q-f_8^4-f_4^2}, respectively.   
\end{proof}

\begin{lemma} We have
\begin{align}
 q \frac{f_4^8}{f_2^4}
&\equiv  2 \frac{f_1^3 f_2^3}{f_3 f_6}
   - 2  \frac{f_1^8}{f_2^4} \pmod{9},   \label{q-f4^8-f2^4} \\
2 \frac{f_2^{20}}{f_1^8 f_4^8}
& \equiv \frac{f_1^8}{f_2^4} + \frac{f_2^{12}f_3^2}{f_1^6 f_6^4}
+ 18 q\frac{f_2^4 f_6^4}{f_1^2 f_3^2} \pmod{27}. \label{2f2^20-f1^8f4^8}
\end{align}
\end{lemma}
\begin{proof}
 Berndt~\cite[p. 40, Entry 25]{Ramanujan-3} also records 
\begin{align}
\varphi^4(q) - \varphi^4(-q) &= 16q \psi^4(q^2). \label{Ram-III-i3}
\end{align}
By \eqref{psi}--\eqref{varphi-q}, \eqref{Ram-III-i3} becomes
\begin{align}\label{Ram-III-i3-f}
\frac{f_2^{20}}{f_1^8 f_4^8} - \frac{f_1^8}{f_2^4} = 16q \frac{f_4^8}{f_2^4}.
\end{align}
We also use Berndt's identity~\cite[p.~150, Eq.~(9.14)]{Ramanujan-5}: 
\begin{align}\label{Ram-id-1}
\frac{\psi^6(q)}{\psi^2(q^3)} + 18q\psi^2(q)\psi^2(q^3)
- 27q^2\frac{\psi^6(q^3)}{\psi^2(q)}
= \varphi^4 (q)+ 16q\psi^4(q^2).
\end{align}
Putting \eqref{psi}--\eqref{varphi-q} into \eqref{Ram-id-1}, then
reducing modulo 9 and modulo 27, respectively, and applying Lemma \ref{f-3-r}, we deduce that
\begin{align}
\frac{f_1^3 f_2^3}{f_3 f_6}
& \equiv  \frac{f_2^{20}}{f_1^8 f_4^8} + 16 q \frac{f_4^8}{f_2^4}   \pmod{9}, \label{Ram-id-1-f-9}  \\
\frac{f_2^{12}f_3^2}{f_1^6 f_6^4} + 18 q \frac{f_2^4 f_6^4}{f_1^2 f_3^2}
& \equiv  \frac{f_2^{20}}{f_1^8 f_4^8} + 16 q \frac{f_4^8}{f_2^4}   \pmod{27}. \label{Ram-id-1-f-27}
\end{align}
Adding \eqref{Ram-III-i3-f} and \eqref{Ram-id-1-f-9}, multiplying both sides by 2, and reducing modulo 9 gives \eqref{q-f4^8-f2^4}.
Subtracting \eqref{Ram-III-i3-f} from \eqref{Ram-id-1-f-27}  yields \eqref{2f2^20-f1^8f4^8}.    \end{proof}

\begin{lemma} We have
\begin{align}
& - 72 \frac{f_2^8 f_3 f_6^4}{f_1^7 f_4^9}
+ 54 q \frac{f_1^{15} f_4^7}{f_2^5f_3f_6} \nonumber \\
& \equiv
- 36 \frac{f_3^3}{f_1^5 f_4}
- 36 \frac{f_1^9 f_3 f_6^4}{f_2^{16}f_4}
+ 81 q \frac{f_2^{16}}{f_1^4 f_4}
+ 108 \frac{f_3^{4}f_6}{f_2^7 f_4}
- 108 \frac{f_3^{8}}{f_1^4f_2^5 f_4 f_6}
\pmod{243}.  \label{gf-con-rel-4-t3}
\end{align}
\end{lemma}
\begin{proof}
 Applying \eqref{2f2^20-f1^8f4^8}, \eqref{q-f4^8-f2^4}, and Lemma \ref{f-3-r}, we obtain
\begin{align*}
& - 72 \frac{f_2^8 f_3 f_6^4}{f_1^7 f_4^9}
+ 54 q \frac{f_1^{15} f_4^7}{f_2^5f_3f_6}   \nonumber \\
& \equiv
- 36 \frac{f_1 f_3 f_6^4}{f_2^{12}f_4} \left(2\frac{f_2^{20}}{f_1^8 f_4^8}\right)
+ 54 \frac{f_1^{15}}{f_2 f_3 f_4 f_6} \left(q \frac{f_4^8}{f_2^4}\right) \nonumber \\
& \equiv
- 36 \frac{f_3^3}{f_1^5 f_4}
- 36 \frac{f_1^9 f_3 f_6^4}{f_2^{16}f_4}
+ 81 q \frac{f_6^8}{f_1 f_2^8 f_3 f_4}
+ 108 \frac{f_1^{18}f_2^2}{f_3^2 f_4 f_6^2}
- 108 \frac{f_1^{23}}{f_2^5 f_3 f_4 f_6}  \nonumber \\
& \equiv
- 36 \frac{f_3^3}{f_1^5 f_4}
- 36 \frac{f_1^9 f_3 f_6^4}{f_2^{16}f_4}
+ 81 q \frac{f_2^{16}}{f_1^4 f_4}
+ 108 \frac{f_3^{4}f_6}{f_2^7 f_4}
- 108 \frac{f_3^{8}}{f_1^4f_2^5 f_4 f_6}
\pmod{243}.
\end{align*}
\end{proof}

\begin{lemma} We have
\begin{align}
  1 - \frac{f_1^{16}}{f_2^8} - q \frac{f_2^{16}}{f_1^8}
& \equiv  0 \pmod{3}, \label{q-f2^{16}-f1^8}  \\
   \frac{f_2f_3^6}{f_1^2f_6^3}
 - \frac{f_3f_6}{f_1^3f_2^3}
 + q\frac{f_1f_6^6}{f_2^2f_3^3}
& \equiv  0  \pmod{9},  \label{f2f3^6-f1^2f6^3}  \\
36 \frac{f_3^{4}f_6}{f_1^4 f_2^7}
- 36 \frac{f_1^5 f_3 f_6^4}{f_2^{16}}
- 81 q f_2^{12}
& \equiv  0 \pmod{243}, \label{81q-f2^{12}} \\
 72 \frac{f_1^{10}f_2^{14}}{f_3^2 f_4^5 f_6^2}
- 72 \frac{f_1 f_3 f_6^4}{f_2^4 f_4^5}
- 81 q \frac{f_2^{24}}{f_1^4 f_4^5}
& \equiv  0 \pmod{243}. \label{gf-con-rel-4-t2}
\end{align}
\end{lemma}
\begin{proof}
We reduce \eqref{q-f1-f6-6} modulo 3 and modulo 9, then apply Lemma \ref{f-3-r} to obtain \eqref{q-f2^{16}-f1^8} and \eqref{f2f3^6-f1^2f6^3}, respectively.
By \eqref{lem-f1-3-f6-2} and Lemma \ref{f-3-r},
\begin{align*}
36 \frac{f_3^{4}f_6}{f_1^4 f_2^7}
- 36 \frac{f_1^5 f_3 f_6^4}{f_2^{16}}
- 81 q f_2^{12}
& \equiv 36 \frac{f_1^2 f_3^2 f_6^2}{f_2^{10}}
\left(\frac{f_2^3 f_3^2}{f_1^6 f_6}
- \frac{f_1^3 f_6^2}{f_2^6 f_3} \right) - 81 q f_2^{12} \nonumber \\
& \equiv
81 q \frac{f_2^6 f_3^2 f_6^2}{f_1^6} - 81 q f_2^{12} \nonumber \\
& \equiv  0 \pmod{243}.
\end{align*}
Thus, \eqref{81q-f2^{12}} follows.
Using \eqref{lem-f1-6-f6-4} and Lemma \ref{f-3-r} gives 
\begin{align*}
 72 \frac{f_1^{10}f_2^{14}}{f_3^2 f_4^5 f_6^2}
- 72 \frac{f_1 f_3 f_6^4}{f_2^4 f_4^5}
- 81 q \frac{f_2^{24}}{f_1^4 f_4^5}
& \equiv
72 \frac{f_1^{10}f_2^{14}}{f_3^2 f_4^5 f_6^2}
- 72 \frac{f_2^8 f_3^3}{f_1^5 f_4^5}
  \left(\frac{f_1^6 f_6^4}{f_2^{12}f_3^2}\right)
- 81 q \frac{f_2^{24}}{f_1^4 f_4^5}   \nonumber \\
& \equiv
 72 \frac{f_1^{10}f_2^{14}}{f_3^2 f_4^5 f_6^2}
- 72 \frac{f_2^{14}f_3^7}{f_1^{17}f_4^5 f_6^2} 
+  81 q \frac{f_2^{27}f_3^5}{f_1^{19}f_4^5f_6}
- 81 q \frac{f_2^{24}}{f_1^4 f_4^5}    \nonumber \\
& \equiv  0 \pmod{243}.
\end{align*}
Hence \eqref{gf-con-rel-4-t2} holds.    
\end{proof}

\section{The first initial generating function}
For a formal series $F(q)=\sum_{n\ge0}a_nq^n$, define
\[
 E_\varepsilon F(q)
   :=\sum_{n\ge0}a_{2n+\varepsilon}q^n
   \qquad(\varepsilon\in\{0,1\}).
\]
Thus $E_0$ extracts the even coefficients, while
$E_1$ extracts the odd coefficients; in each case the
remaining powers of $q^2$ are replaced by powers of $q$.
We use this notation for the successive two-dissections below.

\begin{proof}[Proof of \eqref{1-gf-96}] 
We begin with the generating function of $\overline{p}(3n+1)$ given by Hirschhorn and Sellers \cite{1-Hirschhorn-Sellers-2005}, namely,
\begin{align}\label{gf-3}
\sum_{n=0}^{\infty}\overline{p}(3n+1)q^n
= 2f_2^3(f_1^2)\Big(\frac{f_3}{f_1^3}\Big)^3.
\end{align}
Substituting \eqref{f1-2} and \eqref{f3-f1-3} into \eqref{gf-3}, applying $E_1$, and using Lemma~\ref{f-3-r}, we obtain
\begin{align}\label{gf-6n+4}
&\sum_{n=0}^{\infty}\overline{p}(6n+4)q^n \nonumber \\
&\equiv  18\frac{f_2^{12}f_3^7f_4^5}{f_1^{21}f_6^2f_8^2}
- 4\frac{f_2^{18}f_3^9f_8^2}{f_1^{23}f_4f_6^6}
+ 54q\frac{f_2^{4}f_3^3f_4^5f_6^6}{f_1^{17}f_8^2}
- 108q\frac{f_2^{10}f_3^5f_6^2f_8^2}{f_1^{19}f_4} \nonumber \\
&\equiv  18\frac{f_2^{12}f_4^5}{f_6^2f_8^2}\Big(\frac{f_1^3}{f_3}\Big)^2
 - 4\frac{f_2^{18}f_8^2}{f_4f_6^6}
     \Big(\frac{f_3}{f_1^3}\Big)^9\big(f_1^4\big)
 + 54q\frac{f_2^{4}f_4^5f_6^6}{f_8^2}\frac{f_1}{f_3^3}
 - 108q\frac{f_2f_6^5f_8^2}{f_4}
     \frac{f_1^3}{f_3}\frac{1}{f_1^4} \pmod{243}.
\end{align}
Plugging \eqref{f1-4}, \eqref{1-f1-4}, \eqref{f3-f1-3}, \eqref{f1-3-f3}, and \eqref{f1-f3-3} into \eqref{gf-6n+4},  then applying $E_0$ and Lemma \ref{f-3-r}, we get 
\begin{align}\label{gf-12n+4}
&\sum_{n=0}^{\infty}\overline{p}(12n+4)q^n  \nonumber \\
&\equiv -4\frac{f_2^{63}f_3^{21}}{f_1^{65}f_4^{2}f_6^{18}}
  + 18\frac{f_1^{12}f_2^{11}}{f_3^2f_4^2f_6^2}
  + 54q\frac{f_2^{4}f_3^{5}f_4^6}{f_1^9f_6}
  - 54q\frac{f_1^{7}f_2^{3}f_6^6}{f_3^3f_4^2}
  - 54q\frac{f_2^{47}f_3^{19}f_4^6}{f_1^{59}f_6^{14}}
  \nonumber \\
&~~~
  - 81q\frac{f_2^{55}f_3^{17}}{f_1^{61}f_4^2f_6^{10}}
  + 81q\left\{\frac{f_2^{12}f_3^3f_6^3}{f_1^{11}f_{4}^2}
  - \frac{f_1^{16}f_2^{3}f_6^6}{f_3^{6}f_{4}^2}\right\}
  + 81q^2\frac{f_2^{39}f_3^{15}f_4^6}{f_1^{55}f_6^{6}}  \nonumber \\
&\equiv
 - 4\frac{f_2^{63}}{f_4^{2}f_6^{18}}
         \Big(\frac{f_3}{f_1^3}\Big)^{21}\Big(\frac{1}{f_1^2}\Big)
 + 18\frac{f_2^{11}}{f_4^2f_6^2}\Big(\frac{f_1^3}{f_3}\Big)^3
     \Big(\frac{f_3}{f_1}\Big) \big(f_1^4\big)
 + 54q\frac{f_2^{4}f_4^6}{f_6}\Big(\frac{f_3}{f_1^3}\Big)^3
        \Big(\frac{f_3^3}{f_1}\Big)\Big(\frac{f_1}{f_3}\Big)
 \nonumber \\
&~~~  - 54q\frac{f_2^{3}f_6^6}{f_4^2}\frac{1}{f_1^2}
 - 54qf_2^{2}f_4^{6}f_6\frac{f_3}{f_1^3}\frac{1}{f_1^2}
 - 81q\frac{f_2f_6^{8}}{f_4^2}
  \frac{f_1}{f_3^3}\frac{1}{f_1^2}
 + 81q^2f_4^{6}f_6^{7}\frac{f_1}{f_3^3}\frac{1}{f_1^2}
   \pmod{243}.
\end{align}
Substituting \eqref{f1-4}, \eqref{1-f1-2}, and \eqref{f3-f1-3}--\eqref{f3-f1} into \eqref{gf-12n+4}, we apply $E_0$ and Lemma \ref{f-3-r} to derive 
\begin{align}\label{gf-24n+4}
&\sum_{n=0}^{\infty}\overline{p}(24n+4)q^n \nonumber \\
&\equiv -4\frac{f_2^{124}f_3^{45}f_4^{5}}{f_1^{131}f_6^{42}f_8^2}
   + 18\frac{f_1^{7}f_2^{18}f_8f_{12}^{2}}{f_3f_4^5f_6^{4}f_{24}}
   - 27q\frac{f_2^{116}f_3^{41}f_4^5}{f_1^{127}f_6^{34}f_8^2}
   - 81q\frac{f_1^{13}f_2^2f_4^3f_8f_{12}^2}{f_3^{3}f_{24}}
   - 108q\frac{f_3^6f_8^2}{f_1^{2}f_4}
    \nonumber \\
&~~~
  -72q\frac{f_1^{11}f_2^5f_4^6f_{24}}{f_3f_6^{3}f_8f_{12}}
  -54q\frac{f_2^{26}f_3^8f_4^2f_{24}}{f_1^{24}f_6^{6}f_8f_{12}}
 + \left\{81q\frac{f_2^2f_3f_6^2f_8^2}{f_1^{3}f_4}
    - 18q\frac{f_2^{122}f_3^{43}f_8^2}{f_1^{129}f_{4}f_6^{38}}
    - 108q\frac{f_2^{14}f_3^{4}f_8^2}{f_1^{12}f_{4}f_6^{2}}
   \right\}         \nonumber \\
&~~~ + \left\{ 81q\frac{f_2^{8}f_3^2f_4^5f_6^2}{f_1^{10}f_{8}^2}
   + 81q\frac{f_4^{5}f_6^6}{f_1f_2^4f_3f_{8}^2}
   + 81q\frac{f_2^{8}f_4^{5}f_6^2}{f_1^{4}f_{8}^2}
   \right\}
   + 54q\frac{f_2^{23}f_3^6f_8f_{12}^2}{f_1^{22}f_4f_6^{3}f_{24}}
   + 81q\frac{f_1^{9}f_2^{13}f_6f_{24}}{f_3^{3}f_{4}^2f_8f_{12}}
   \nonumber \\
&~~~
    + \left\{81q^2\frac{f_2^{6}f_6^{6}f_8^2}{f_1^{2}f_{3}^2f_4}
    - 54q^2\frac{f_2^{114}f_3^{39}f_8^2}{f_1^{125}f_4f_6^{30}}
   \right\}   \nonumber \\
&\equiv -4\frac{f_2^{124}f_4^{5}}{f_6^{42}f_8^2}
  \Big(\frac{f_3}{f_1^3}\Big)^{45} \big(f_1^4\big)
   + 18\frac{f_6^{5}f_8f_{12}^{2}}{f_2^9f_4^5f_{24}}
   \frac{f_1^3}{f_3} f_1^4
   - 27q\frac{f_4^{5}f_6^{5}}{f_2f_8^2}
    \Big(\frac{f_3}{f_1^3}\Big)^2 \Big(\frac{1}{f_1^4}\Big)
   - 81q\frac{f_2^2f_4^9}{f_8^2} f_1^4
     \nonumber \\
&~~~- 108q\frac{f_8^2}{f_4}\Big(\frac{f_3^3}{f_1}\Big)^2
    - 72q\frac{f_2^5f_4^6f_{24}}{f_6^{3}f_8f_{12}}
    \Big(\frac{f_1^3}{f_3}\Big)^4 \Big(\frac{f_3^3}{f_1}\Big)
    -54q\frac{f_4^2f_6^3f_{24}}{f_2f_8f_{12}}\frac{f_1^3}{f_3}
    -45q\frac{f_2^{14}f_8^2}{f_4f_6^2}\Big(\frac{f_1^3}{f_3}\Big)^2
        \nonumber \\
&~~~ + 54q\frac{f_2^{5}f_6^3f_8f_{12}^2}{f_4f_{24}}\frac{1}{f_1^4}
+ 81q\frac{f_2f_6^{5}f_{8}^2}{f_{4}^5}
    + 27q^2\frac{f_6^{9}f_8^2}{f_2^{3}f_4}\frac{f_1}{f_3^3}
   \pmod{243}.
\end{align}
We use \eqref{f1-4}, \eqref{1-f1-4}, and \eqref{f3-f1-3}--\eqref{f1-f3-3} in
\eqref{gf-24n+4}, then apply $E_1$ and Lemma~\ref{f-3-r} to obtain 
\begin{align}\label{gf-48n+28}
&\sum_{n=0}^{\infty}\overline{p}(48n+28)q^n \nonumber \\
&\equiv 16\frac{f_2^{273}f_3^{93}f_{4}^2}{f_{1}^{279}f_6^{90}}
  - 54\frac{f_2^{281}f_3^{91}}{f_{1}^{281}f_{4}^6f_6^{86}}
  - 27\frac{f_2^{31}f_3^{11}}{f_1^{33}f_4^6f_6^4}
  - 54\frac{f_2^5f_3^3f_{12}}{f_1f_4f_6^2}
  - 72\frac{f_3^5f_4^5f_6}{f_1^7f_2^4f_{12}}
  - 72\frac{f_1^3f_2^{21}f_{12}}{f_3f_4f_6^6}
  \nonumber \\
&~~~ - 81\frac{f_2^{19}}{f_4^6}
  + 81\frac{f_1f_3^5f_4^2}{f_2^5}
  - \left\{108\frac{f_2^5f_3^4f_4^2}{f_1^4f_6^2}
  + 45\frac{f_1^{14}f_2^{5}f_4^2}{f_3^2f_6^2}
  \right\}
  +  \left\{54 \frac{f_2^{13}f_3^3f_6^2}{f_1^9f_4^3f_{12}}
  - 54\frac{f_2^4f_3^3f_6^5}{f_1^9f_4^3f_{12}}
  \right\} \nonumber \\
&~~~
  - \left\{ 81q\frac{f_2^{265}f_3^{89}f_4^2}{f_1^{275}f_6^{82}}
     + 108q\frac{f_1^7f_2^{13}f_6^2f_{12}}{f_3^5f_4}
     \right\}
  + 81q\frac{f_2^{273}f_3^{87}}{f_1^{277}f_4^6f_6^{78}}
  \nonumber \\
&~~~
  + \left\{81q\frac{f_1^{18}f_4^{2}f_6^6}{f_2^3f_3^6}
  + 81q\frac{f_2^{15}f_3^9f_4^2}{f_1^{27}}
  + 108q\frac{f_4^{2}f_6^6}{f_2^{3}}
  \right\}  \nonumber \\
&\equiv 16\frac{f_2^{30}f_{4}^2}{f_6^{9}}\Big(\frac{f_3}{f_1^3}\Big)^{12}
  - 54\frac{f_2^{2}f_6^{7}}{f_{4}^6}
\frac{f_1}{f_3^3}\frac{f_3}{f_1^3}
  - 27\frac{f_2^{4}f_6^{5}}{f_4^6}\frac{f_1^3}{f_3}
 - 54\frac{f_2^5f_{12}}{f_4f_6^2}\frac{f_3^3}{f_1}
  - 72\frac{f_4^5f_6}{f_2^4f_{12}}
    \Big(\frac{f_3}{f_1^3}\Big)^2\Big(\frac{f_3^3}{f_1}\Big)
     \nonumber \\
&~~ - 72\frac{f_6^{3}f_{12}}{f_2^6f_{4}}\frac{f_1^3}{f_3}
  - 81\frac{f_2f_6^{6}}{f_{12}^2}
  + 81\frac{f_4^2}{f_2^2f_6}\Big(\frac{f_3^3}{f_1}\Big)^2
  + 90\frac{f_2^{5}f_4^{2}}{f_6^2}
        \Big(\frac{f_1^3}{f_3}\Big)^5\Big(\frac{f_3^3}{f_1}\Big)
    \nonumber \\
&~~+ 54q\frac{f_2^{4}f_6^{5}f_{12}}{f_4}
    \Big(\frac{f_1}{f_3^3}\Big) \Big(\frac{f_1^3}{f_3}\Big)^2
  + 81q\frac{f_6^{13}}{f_{12}^{2}}\Big(\frac{f_1}{f_3^3}\Big)^2
  + 27q\frac{f_4^{2}f_6^6}{f_2^{3}}
   \pmod{243}.
\end{align}
Using \eqref{f3-f1-3}--\eqref{f1-f3-3} in \eqref{gf-48n+28}, we apply $E_1$
and Lemma~\ref{f-3-r} to obtain 
\begin{align}\label{gf-96n+76}
&\sum_{n=0}^{\infty}\overline{p}(96n+76)q^n  \nonumber \\
&\equiv
 \left\{27\frac{f_2^2f_3^6}{f_1^3}
  + 81\frac{f_1^6f_3^3f_6^3}{f_2^7} \right\}
+ \left\{90\frac{f_2^{70}f_3^{25}}{f_1^{76}f_6^{20}}
  - 45\frac{f_1^5f_2^{16}}{f_3^2f_6^2}
  - 18\frac{f_2^{16}f_3^7}{f_1^{22}f_6^2}
  \right\} \nonumber \\
&~~~ + 108\frac{f_3f_6^{4}}{f_1^4f_2^2}
+ \left\{54\frac{f_1^5f_2^7f_6}{f_3^2}
  - 54\frac{f_1^5f_6^4}{f_2^2f_3^2} \right\}
+ \left\{81\frac{f_1^2f_2^4f_6^2}{f_3}
  - 81\frac{f_2^4f_3f_6^2}{f_1^4} \right\}   \nonumber \\
&~~~ + \left\{27q\frac{f_2^{62}f_3^{21}}{f_1^{72}f_6^{12}}
  + 81q\frac{f_1^6f_6^{10}}{f_2^4f_3^5}
  + 81q\frac{f_1^9f_2^8f_6^{6}}{f_3^6}
  + 81q\frac{f_2^8f_3^3f_6^{6}}{f_1^{18}}
  + 81q\frac{f_1^9f_6^{9}}{f_2f_3^6}
  \right\}  \nonumber \\
&\equiv
108 \frac{f_2^2f_3^6}{f_1^3}
  + 27\frac{f_1^5f_2^{16}}{f_3^2f_6^2}
  + 108\frac{f_3f_6^{4}}{f_1^4f_2^2}
  + 108q\frac{f_6^9}{f_2f_3^3}     \pmod{243}.
\end{align}
By Lemma \ref{f-3-r}, it follows that
\begin{align}\label{gf-96n-coeff-sim}
\frac{f_1^5f_2^{16}}{f_3^2f_6^2}
\equiv\frac{f_3f_6^{4}}{f_1^4f_2^2} \pmod{9}.
\end{align}
Substituting \eqref{gf-96n-coeff-sim} into \eqref{gf-96n+76} and using $135 \equiv -108 \pmod{243}$, we conclude that
\begin{align}\label{gf-96n-res}
\sum_{n=0}^{\infty}\overline{p}(96n+76)q^n
\equiv
   108 \left( \frac{f_2^2f_3^6}{f_1^3}
 - \frac{f_3f_6^{4}}{f_1^4f_2^2}
 + q\frac{f_6^9}{f_2f_3^3} \right)  \pmod{243}.
\end{align}
\end{proof}

\section{The second initial generating function}

\begin{proof}[Proof of \eqref{2-gf-96}] 
Applying $E_1$ to \eqref{gf-12n+4} and using Lemma~\ref{f-3-r} gives 
\begin{align}\label{gf-24n+16}
&\sum_{n=0}^{\infty}\overline{p}(24n+16)q^n  \nonumber \\
&\equiv
- 8\frac{f_2^{126}f_3^{45}f_8^{2}}{f_1^{131}f_4 f_6^{42}}
- 9\frac{f_2^{120}f_3^{43}f_4^{5}}{f_1^{129}f_6^{38}f_8^2}
+ 18\frac{f_1^{7}f_2^{17}f_{24}}{f_3f_4^2f_6^{3}f_8f_{12}}
- 72\frac{f_1^{11}f_2^{6}f_4^{3}f_8f_{12}^2}{f_3f_6^{4}f_{24}}  + 81\frac{f_1^{9}f_2^{14}f_8f_{12}^2}{f_3^3f_4^{5}f_{24}}
  \nonumber \\
&~~
    + 54\frac{f_2^{27}f_3^{8}f_8f_{12}^2}{f_1^{24}f_4f_6^{7}f_{24}}
     - \left\{ 81\frac{f_3f_4^5f_6^2}{f_1^3f_8^2}
     + 54\frac{f_2^{12}f_3^4f_4^5}{f_1^{12}f_6^2f_8^2}
     \right\}
    - 54\frac{f_3^{6}f_4^{5}}{f_1^{2}f_2^{2}f_8^2}
    - 81q\frac{f_1^{13}f_2f_4^{6}f_6f_{24}}{f_3^{3}f_8f_{12}} \nonumber \\
&~~
  - \left\{ 27q\frac{f_2^{112}f_3^{39}f_4^5}{f_1^{125}f_6^{30}f_8^2}
  + 81q\frac{f_2^{4}f_4^5f_6^6}{f_1^{2}f_3^2f_8^2}
     \right\}
  - \left\{ 54q\frac{f_2^{118}f_3^{41}f_8^2}{f_1^{127}f_4f_6^{34}}
    + 81q\frac{f_6^{6}f_8^2}{f_1f_2^{2}f_3f_4}
     \right\} \nonumber \\
&~~
  - \left\{ 54q\frac{f_2^{22}f_3^{6}f_4^2f_{24}}{f_1^{22}f_6^2f_8f_{12}}
    + 81q\frac{f_2^{10}f_3^2f_6^2f_8^2}{f_1^{10}f_4}
    + 81q\frac{f_2^{10}f_6^2f_8^2}{f_1^{4}f_4}
     \right\}   \nonumber \\
&\equiv
 - 8\frac{f_2^{126}f_8^{2}}{f_4 f_6^{42}}
      \Big(\frac{f_3}{f_1^3}\Big)^{45}\big(f_1^4\big)
 - 9\frac{f_2^{12}f_4^{5}}{f_6^{2}f_8^2}\Big(\frac{f_1^3}{f_3}\Big)^{2}
 + 18\frac{f_2^{17}f_{24}}{f_4^2f_6^{3}f_8f_{12}}      \frac{f_1^3}{f_3}f_1^4
 - 72\frac{f_2^{6}f_4^{3}f_8f_{12}^2}{f_6^{4}f_{24}}
       \Big(\frac{f_1^3}{f_3}\Big)^4\Big(\frac{f_3^3}{f_1}\Big)
      \nonumber \\
&~~ + 81\frac{f_4f_6^{5}}{f_2f_8^{2}}
    + 54\frac{f_6^{2}f_8f_{12}^2}{f_4f_{24}}\frac{f_1^3}{f_3}
  + 108\frac{f_2^3f_4^5f_6}{f_8^2}\frac{f_3}{f_1^3}
  - 54\frac{f_4^{5}}{f_2^{2}f_8^2}\Big(\frac{f_3^3}{f_1}\Big)^2
    - 81qf_2^4f_4^{3}f_8^2(f_1^4)    \nonumber \\
&~~ - 108q\frac{f_2^{4}f_4^5f_6^6}{f_8^2}\frac{f_1}{f_3^3}
    + 108q\frac{f_2f_6^{5}f_8^2}{f_4}\frac{1}{f_1f_3}
    + 27q\frac{f_2^{4}f_4^2f_6^4f_{24}}{f_8f_{12}}\frac{1}{f_1^4}
  \pmod{243}.
\end{align}
Using \eqref{f1-4}, \eqref{f3-f1-3}--\eqref{f1-f3-3},
and \eqref{1-f1f3} in \eqref{gf-24n+16}, 
we apply $E_0$ and Lemma~\ref{f-3-r} to get
\begin{align}\label{gf-48n+16}
&\sum_{n=0}^{\infty}\overline{p}(16(3n+1))q^n \nonumber \\
&\equiv - 8\frac{f_2^{279}f_3^{93}}{f_1^{281}f_4^2 f_6^{90}}
     - 9\frac{f_1^{12}f_2^{11}}{f_3^{2}f_4^{2}f_6^2}
     + 18\frac{f_1^{15}f_2^{11}f_{12}}{f_3^3f_4^5f_6^{2}}
     + 54\frac{f_2^2f_3^{2}f_4f_6}{f_{12}}
     - 72\frac{f_1^{4}f_2^{18}f_4}{f_3^2f_6^{3}f_{12}}
     + 54\frac{f_2^{11}f_3^{4}}{f_1^{6}f_4^2f_6^{2}} \nonumber \\
&~~~
     + 81\frac{f_2f_3^5}{f_1f_4^2}
     + \left\{81qf_1^6f_2f_4^6
     + 108q\frac{f_2^{4}f_3^3f_{12}^2}{f_1^{3}f_6}
     + 108q\frac{f_2^{4}f_3^{4}f_4^3f_{12}}{f_1^{6}f_6}
     - 27q\frac{f_1^{21}f_4^3f_6^2f_{12}}{f_2^{5}f_3^5}
       \right\}      \nonumber \\
&~~~ + \left\{108q\frac{f_1^{7}f_2^{3}f_6^6}{f_3^{3}f_4^2}
    - 54q\frac{f_2^{3}f_6^6}{f_1^{2}f_4^2}
    - 81q\frac{f_1^{16}f_2^3f_6^6}{f_3^{6}f_4^2}
     \right\}
    - 54q\frac{f_2^{263}f_3^{91}f_4^6}{f_1^{275}f_6^{86}}
    \nonumber \\
&~~~
    - \left\{108q\frac{f_1^{8}f_2^{10}f_4f_6^5}{f_3^{6}f_{12}}
    + 81q\frac{f_2^{271}f_3^{89}}{f_1^{277}f_4^{2}f_6^{82}}
     \right\}
    + 81q^2\frac{f_2^{255}f_3^{87}f_4^6}{f_1^{271}f_6^{78}}      \nonumber \\
&\equiv - 8\frac{f_2^{36}}{f_4^2 f_6^{9}}
      \Big(\frac{f_3}{f_1^3}\Big)^{12}\Big(\frac{1}{f_1^2}\Big)
      - 9\frac{f_2^{11}}{f_4^{2}f_6^2}
      \Big(\frac{f_1^3}{f_3}\Big)^{4}\big(f_3^2\big)
      + 18\frac{f_2^{11}f_{12}}{f_4^5f_6^{2}}
      \Big(\frac{f_1^3}{f_3}\Big)^{5}\big(f_3^2\big)
      + 54\frac{f_2^2f_4f_6}{f_{12}}f_3^2
       \nonumber \\
&~~~ - 72\frac{f_2^{18}f_4}{f_6^{3}f_{12}}
      \Big(\frac{f_1^3}{f_3}\Big)^{2}\Big(\frac{1}{f_1^2}\Big)
     + 54\frac{f_2^{2}f_6}{f_4^2}
       \frac{f_1^3}{f_3}f_3^2
     + 81\frac{f_2}{f_4^2}\frac{f_3^3}{f_1}f_3^2
     + 27q\frac{f_2^{4}f_{12}^2}{f_6}
        \frac{f_3^3}{f_1}\frac{1}{f_1^2} \nonumber \\
&~~~  - 27q\frac{f_2^{3}f_6^6}{f_4^2}\frac{1}{f_1^2}
     - 54qf_2^{2}f_4^{6}f_6\frac{f_3}{f_1^3}\frac{1}{f_1^2}
     + 54q\frac{f_2f_4f_6^8}{f_{12}}
       \frac{f_1}{f_3^3}\frac{1}{f_1^2}
     + 81q^2f_4^6f_6^{7}\frac{f_1}{f_3^3}\frac{1}{f_1^2}     \pmod{243}.
\end{align}
Substituting \eqref{1-f1-2} and \eqref{f3-2}--\eqref{f1-f3-3} into
\eqref{gf-48n+16}, applying $E_1$ and  Lemma~\ref{f-3-r}, we obtain
\begin{align}\label{gf-16(6n+4)}
&\sum_{n=0}^{\infty}\overline{p}(16(6n+4))q^n  \nonumber \\
& \equiv
     99\frac{f_1^{13}f_2^9f_8^2}{f_3^3f_4f_6^3}
    - 16\frac{f_2^{72}f_3^{27}f_8^2}{f_1^{77}f_4f_6^{24}}
    - 27\frac{f_3^6f_4^5}{f_1^2f_2^2f_8^2}
+81\left\{\frac{f_1^{13}f_2^{6}f_{12}^5}{f_3^{3}f_6^{2}f_{24}^2}
    + \frac{f_1^{4}f_6f_{12}^5}{f_2^{3}f_{24}^2}
    + \frac{f_1f_3f_6f_{12}^5}{f_2^3f_{24}^2}
       \right\}
    \nonumber \\
&~~
   + \left\{ 81\frac{f_2^3f_3f_4^5f_6}{f_1^3f_8^2}
       - 54 \frac{ f_2^{12}f_3^4f_4^5}{f_1^{12}f_6^2f_8^2}
        - 54 \frac{f_1^{15}f_2^3f_4^5f_6}{f_3^5f_8^2}
    - 45\frac{f_2^{66}f_3^{25}f_4^5}{f_1^{75}f_6^{20}f_8^2}
       \right\}
    - 81q^2\frac{f_1^{15}f_2^2f_6^4f_{24}^2}{f_3^5f_{12}}     \nonumber \\
&~~ + \left\{27q\frac{f_1^{2}f_2f_3^{2}f_{24}^2}{f_6f_{12}}
    - 18q\frac{f_1^{11}f_2^{10}f_{24}^2}{f_3f_6^{4}f_{12}}
    + 81q\frac{f_2f_3^{3}f_{24}^2}{f_1f_6f_{12}}
      \right\}
+ \left\{108q\frac{f_2^{58}f_3^{21}f_4^{5}}{f_1^{71}f_6^{12}f_8^2}
    - 81q\frac{f_2^4f_4^5f_6^{6}}{f_1^2f_3^2f_8^{2}}
      \right\} \nonumber \\
&~~  - 27q\left\{\frac{f_2^{64}f_3^{23}f_8^{2}}{f_1^{73}f_4f_6^{16}}
    + 2\frac{f_2f_6^{5}f_8^{2}}{f_1f_3f_4}
    + 3 \frac{f_2^{10}f_6^{2}f_{8}^2}{f_1^4f_4}
    + 3 \frac{f_1^{17}f_2f_6^{5}f_8^{2}}{f_3^{7}f_4}
    + 3 \frac{f_2^{10}f_3^2f_6^{2}f_{8}^2}{f_1^{10}f_4}
      \right\}
   \nonumber \\
& \equiv
 99\frac{f_2^9f_8^2}{f_4f_6^3}
      \Big(\frac{f_1^3}{f_3}\Big)^{3}\big(f_1^4\big)
 - 16\frac{f_2^{72}f_8^2}{f_4f_6^{24}}
      \Big(\frac{f_3}{f_1^3}\Big)^{27}\big(f_1^4\big)
 - 27\frac{f_4^5}{f_2^2f_8^2}\Big(\frac{f_3^3}{f_1}\Big)^{2}
 - 72\frac{f_2^{12}f_4^5}{f_6^2f_8^2}
      \Big(\frac{f_1^3}{f_3}\Big)^{2} \nonumber \\
&~~ - 81q^2\frac{f_2^2f_6^4f_{24}^2}{f_{12}}
  + 90q\frac{f_2^{10}f_{24}^2}{f_6^{4}f_{12}}
        \Big(\frac{f_1^3}{f_3}\Big)^{4}\Big(\frac{f_3^3}{f_1}\Big)
  + 27q\frac{f_2^4f_4^5f_6^{6}}{f_8^{2}}\frac{f_1}{f_3^3}
  - 81q\frac{f_2f_6^{5}f_8^{2}}{f_4}\frac{1}{f_1^4}
      \pmod{243}.
\end{align}
Using \eqref{f1-4}, \eqref{1-f1-4}, and \eqref{f3-f1-3}--\eqref{f1-f3-3}
in \eqref{gf-16(6n+4)}, we apply $E_0$ and Lemma~\ref{f-3-r} to obtain 
\begin{align}\label{gf-192n+64}
& \sum_{n=0}^{\infty}\overline{p}(16(12n+4))q^n   \nonumber \\
& \equiv
  99\frac{f_1^7f_2^{18}}{f_3^3f_4^2f_6^3}
  - 16\frac{f_2^{171}f_3^{57}}{f_1^{173}f_4^2f_6^{54}}
  - \left\{27\frac{f_2^{11}f_3^4}{f_1^6f_4^2f_6^2}
  + 72\frac{f_1^{12}f_2^{11}}{f_3^2f_4^2f_6^2}  \right\}
  - 81q\frac{f_2f_3^5f_{4}^6}{f_1^9}
  - 18q\frac{f_1^{10}f_2^{11}f_{12}^2}{f_3^4f_6^2}
  \nonumber \\
&~~~
  + 27q\left\{3\frac{f_1^{16}f_2^{3}f_6^6}{f_3^6f_4^2}
  - \frac{f_2^{3}f_6^6}{f_1^2f_4^2}
  - \frac{f_1^7f_2^{3}f_6^6}{f_3^3f_4^2}
   \right\}
   + 81q \left\{\frac{f_2^{155}f_3^{55}f_4^6}{f_1^{167}f_6^{50}}
  - \frac{f_1^{13}f_2^2f_4^6f_6}{f_3^5}  \right\}
  - 81q\frac{f_1^{2}f_3^4f_{12}^2}{f_6}
   \nonumber \\
& \equiv
   99\frac{f_2^{18}}{f_4^2f_6^3}
      \Big(\frac{f_1^3}{f_3}\Big)^{3}\Big(\frac{1}{f_1^2}\Big)
   - 16\frac{f_6^{27}}{f_2^{72}f_4^{2}}
      \Big(\frac{f_1^3}{f_3}\Big)^{24}\Big(\frac{1}{f_1^2}\Big)
   - 99\frac{f_2^{11}}{f_4^2f_6^2}
      \Big(\frac{f_1^3}{f_3}\Big)^{4}\big(f_3^2\big)
   - 81qf_2f_{12}^2\big(f_3^2\big)
  \nonumber \\
&~~~ - 18q\frac{f_2^{11}f_{12}^2}{f_6^2}
      \Big(\frac{f_1^3}{f_3}\Big)^{4}\Big(\frac{1}{f_1^2}\Big)
   + 27q\frac{f_2^{3}f_6^6}{f_4^2}\frac{1}{f_1^2}
   - 81q\frac{f_{12}^2}{f_6}\frac{f_3^3}{f_1}f_3^2
       \pmod{243}.
\end{align}
Substituting \eqref{1-f1-2}, \eqref{f3-2}, \eqref{f1-3-f3}, and \eqref{f3-3-f1} into
\eqref{gf-192n+64}, applying $E_0$, and using Lemma~\ref{f-3-r}, we obtain 
\begin{align}\label{gf-384n+64}
&\sum_{n=0}^{\infty}\overline{p}(16(24n+4))q^n \nonumber \\
& \equiv
  - 16\frac{f_2^{70}f_3^{27}f_4^5}{f_1^{77}f_6^{24}f_8^2}
  + 99\frac{f_1^{13}f_2^7f_4^5}{f_3^3f_6^3f_8^2}
  - 99\frac{f_1^{11}f_2^{10}f_{12}^5}{f_3f_6^6f_{24}^2}
  + 54q\frac{f_3^6f_8^2}{f_1^2f_4}
  - 81q\frac{f_6^3f_{12}^5}{f_2f_{24}^2} \nonumber \\
&~~~ + \left\{117q\frac{f_2^{68}f_3^{25}f_8^2}{f_1^{75}f_4f_6^{20}}
   - 36q\frac{f_1^6f_2^{14}f_{8}^2}{f_3^2f_4f_6^2}
   - 81q\frac{f_1^{15}f_2^5f_6f_8^2}{f_3^5f_4}
      \right\}
+\left\{108q\frac{f_2^{62}f_3^{23}f_4^5}{f_1^{73}f_6^{16}f_8^2}
   - 27q\frac{f_1^8f_2^{8}f_{4}^5f_6^2}{f_3^4f_8^2}
      \right\} \nonumber \\
&~~~
   + 108q^2\frac{f_2^{60}f_3^{21}f_8^2}{f_1^{71}f_4f_6^{12}}
 + \left\{54q^2\frac{f_1^{13}f_2^{6}f_{24}^2}{f_3^{3}f_{12}}
   - 81q^2\frac{f_1f_3f_{6}^2f_{24}^2}{f_{12}}
   - 81q^2\frac{f_2^3f_3^2f_{6}f_{24}^2}{f_1^2f_{12}}
      \right\}   \nonumber \\
& \equiv
   - 16\frac{f_2^{70}f_4^5}{f_6^{24}f_8^2}
      \Big(\frac{f_3}{f_1^3}\Big)^{27}\big(f_1^4\big)
  + 99\frac{f_2^7f_4^5}{f_6^3f_8^2}\Big(\frac{f_1^3}{f_3}\Big)^{3}\big(f_1^4\big)
 - 99\frac{f_2^{10}f_{12}^5}{f_6^6f_{24}^2}
      \Big(\frac{f_1^3}{f_3}\Big)^{4}\Big(\frac{f_3^3}{f_1}\Big)
    + 54q\frac{f_8^2}{f_4}\Big(\frac{f_3^3}{f_1}\Big)^2
  \nonumber \\
&~~~     - 81q\frac{f_6^3f_{12}^5}{f_2f_{24}^2}
+ 81q\frac{f_4^{5}f_{6}^5}{f_2f_8^2}\frac{1}{f_1^4}
     + 108q^2\frac{f_6^{9}f_8^2}{f_2^{3}f_4}\frac{f_1}{f_3^3}
     - 108q^2\frac{f_6^3f_{24}^2}{f_2^3f_{12}}f_1^4
      \pmod{243}.
\end{align}
Using \eqref{f1-4}, \eqref{1-f1-4}, and \eqref{f3-f1-3}--\eqref{f1-f3-3}
in \eqref{gf-384n+64}, we apply $E_1$ and Lemma~\ref{f-3-r} to obtain
\begin{align}\label{gf-768n+448}
&\sum_{n=0}^{\infty}\overline{p}(16(48n+28))q^n \nonumber \\
& \equiv
    64\frac{f_2^{165}f_3^{57}f_4^2}{f_1^{171}f_6^{54}}
  + 117\frac{f_1^{10}f_2^{11}f_6^4}{f_3^6f_{12}^2}
  + 54\frac{f_2^{5}f_3^{4}f_{4}^2}{f_1^4f_6^2}
  + 90\frac{f_1^9f_2^{12}f_4^2}{f_3^3f_6^3}
  + 81\frac{f_2^{19}f_3^5}{f_1^{15}f_4^6}
  \nonumber \\
&~~~
  - 81\frac{f_3^3f_6^{5}}{f_1f_{12}^2}
  + \left\{81\frac{f_1^{7}f_2^{20}f_6}{f_3^{5}f_4^{6}}
  - 81\frac{f_2^{173}f_3^{55}}{f_1^{173}f_4^{6}f_{6}^{50}}
      \right\}
- 54q\frac{f_3^3f_4^4f_{12}^2}{f_1f_2^2f_6}
  - 54q\frac{f_4^2f_6^6}{f_2^3}  \nonumber \\
& \equiv
    64\frac{f_4^2f_6^{27}}{f_2^{78}}\Big(\frac{f_1^3}{f_3}\Big)^{24}
  + 117\frac{f_2^{11}f_6^4}{f_{12}^2}
     \Big(\frac{f_1^3}{f_3}\Big)^{3}\Big(\frac{f_1}{f_3^3}\Big)
  + 54\frac{f_2^{5}f_{4}^2}{f_6^2}
      \frac{f_3^3}{f_1}\frac{f_3}{f_1^3}
  +90\frac{f_2^{12}f_4^2}{f_6^3}\Big(\frac{f_1^3}{f_3}\Big)^{3}
   \nonumber \\
&~~~
  + 81\frac{f_2f_6^6}{f_{12}^2}
  - 81\frac{f_6^{5}}{f_{12}^2}\frac{f_3^3}{f_1}
  -54q\frac{f_4^4f_{12}^2}{f_2^2f_6}\frac{f_3^3}{f_1}
  - 54q\frac{f_4^2f_{6}^6}{f_2^3}  \pmod{243}.
\end{align}
Using \eqref{f3-f1-3}--\eqref{f1-f3-3} in \eqref{gf-768n+448}, we apply $E_1$
and Lemma~\ref{f-3-r} to obtain 
\begin{align}\label{gf-1536n+1216}
&\sum_{n=0}^{\infty}\overline{p}(16(96n+76))q^n \nonumber \\
& \equiv
  - \left\{81\frac{f_3^{5}f_6}{f_2}
     + 54\frac{f_2^{2}f_3^{6}}{f_1^{3}}
      \right\}
  - \left\{81\frac{f_2^{7}f_3f_6}{f_1^4}
    + 36\frac{f_1^{14}f_2^{7}f_6}{f_3^{5}}
      \right\}
    + 9\frac{f_2^{70}f_3^{25}}{f_1^{76}f_6^{20}}
    + 27q\frac{f_2^{62}f_3^{21}}{f_1^{72}f_6^{12}}   \nonumber \\
& \equiv
     108\frac{f_2^{2}f_3^{6}}{f_1^{3}}
   - 117\frac{f_1^{14}f_2^{7}f_6}{f_3^{5}}
   + 9\frac{f_1^{5}f_2^{16}}{f_3^{2}f_6^{2}}
   + 27q\frac{f_6^{9}}{f_2f_3^{3}}  \pmod{243}.
\end{align}
Using \eqref{lem-f1-3-f6-2} and Lemma~\ref{f-3-r}, the remaining terms simplify as follows: 
\begin{align}\label{gf-1536n+1216-ter}
&  - 117 \frac{f_1^{14}f_2^{7}f_6}{f_3^{5}}
  + 9 \frac{f_1^{5} f_2^{16}}{f_3^{2}f_6^{2}}
   + 27q\frac{f_6^{9}}{f_2f_3^{3}} \nonumber \\
& \equiv
 - 117 \frac{f_1^{11}f_2^{13}}{f_3^{4}f_6} \left(\frac{f_1^3 f_6^2}{f_2^6 f_3} \right)
 + 9 \frac{f_1^{5} f_2^{16}}{f_3^{2}f_6^{2}}
 + 27q\frac{f_6^{9}}{f_2f_3^{3}} \nonumber \\
& \equiv
 - 117 \frac{f_1^{11}f_2^{13}}{f_3^{4}f_6}
 \left(\frac{f_2^3f_3^2}{f_1^6f_6} - 9q\frac{f_2^{16}}{f_1^8} \right)
 + 9 \frac{f_1^{5} f_2^{16}}{f_3^{2}f_6^{2}}
 + 27q\frac{f_6^{9}}{f_2f_3^{3}}  \nonumber \\
& \equiv
- 108 \frac{f_1^5 f_2^{16}}{f_3^2 f_6^2}
+ \left\{ 81q \frac{f_1^3 f_2^{29}}{f_3^4 f_6}
+ 27q\frac{f_6^{9}}{f_2f_3^{3}} \right\} \nonumber \\
& \equiv
 - 108\frac{f_3f_6^{4}}{f_1^{4}f_2^{2}}
 + 108q\frac{f_6^{9}}{f_2f_3^{3}} \pmod{243}.
\end{align}
Substituting \eqref{gf-1536n+1216-ter} into \eqref{gf-1536n+1216}, we arrive at
\[
\sum_{n=0}^{\infty}\overline{p}(16(96n+76))q^n
\equiv
108\left(\frac{f_2^2f_3^6}{f_1^3}
- \frac{f_3f_6^{4}}{f_1^{4}f_2^{2}}
+ q\frac{f_6^{9}}{f_2f_3^3} \right) \pmod{243}.
\]
\end{proof}

\section{The scaling relation}
We first establish a generating-function congruence for $\overline{p}(16^2(6n+4))$.

\begin{lemma}\label{thm-gf-16^2(6n+4)}
We have
\begin{align}\label{gf-16^2(6n+4)}
& \sum_{n=0}^{\infty}\overline{p}(16^2(6n+4))q^n \nonumber \\
& \equiv
- 64 \frac{f_2^{72}f_3^{27}f_8^2}{f_1^{77}f_4f_6^{24}}
- 108 \frac{f_3^6 f_4^5}{f_1^2 f_2^2 f_8^2}
- 90 \frac{f_1^{13}f_2^9 f_8^2}{f_3^3 f_4 f_6^3}
+ 45 \frac{f_1^6 f_2^{12}f_4^5}{f_3^2 f_6^2 f_8^2}
+ 108 \frac{f_2^4 f_3^8 f_8^2}{f_1^4 f_4 f_6^4}
\nonumber \\
& ~~~
- 27 q \frac{f_2 f_6^5 f_8^2}{f_1 f_3 f_4}
+ 108 q \frac{f_1 f_2^4 f_4^5 f_6^6}{f_3^3 f_8^2}
- 54 q \frac{f_1^2 f_2 f_3^2 f_{24}^2}{f_6 f_{12}}
- 81 q^2 \frac{f_6^5 f_{24}^2}{f_2f_{12}}
  \pmod{243}.
\end{align}
\end{lemma}

\pf Applying $E_1$ to \eqref{gf-192n+64} and then using Lemma~\ref{f-3-r} gives
\begin{align}\label{gf-384n+256}
& \sum_{n=0}^{\infty}\overline{p}(16(24n+16))q^n \nonumber \\
& \equiv
- 32\frac{f_2^{72} f_3^{27} f_8^2}{f_1^{77} f_4 f_6^{24}}
- 45\frac{f_1^{13} f_2^9 f_8^2}{f_3^3 f_4 f_6^3}
+ 27\frac{f_3^6 f_4^5}{f_1^2 f_2^2 f_8^2}
- \left\{
  27\frac{f_1^{13} f_2^6 f_{12}^5}{f_3^3 f_6^2 f_{24}^2}
+ 81\frac{f_1 f_3 f_{12}^5}{f_{24}^2}
+ 81\frac{f_2^3 f_3^2 f_{12}^5}{f_1^2 f_6 f_{24}^2} \right\}  \nonumber \\
& ~~~
+ \left\{
 81\frac{f_1^{15} f_2^3 f_4^5 f_6}{f_3^5 f_8^2}
- 18\frac{f_1^6 f_2^{12}f_4^5}{f_3^2 f_6^2 f_8^2}
- 63\frac{f_2^{66}f_3^{25}f_4^5}{f_1^{75}f_6^{20}f_8^2}  \right\}
- \left\{
 27 q\frac{f_2^{64}f_3^{23}f_8^2}{f_1^{73}f_4f_6^{16}}
+ 54 q\frac{f_1^8 f_2^{10} f_6^2 f_8^2}{f_3^4 f_4}  \right\}
\nonumber \\
& ~~~ - 45 q\frac{f_1^{11}f_2^{10}f_{24}^2}{f_3 f_6^4 f_{12}}
+ 54 q\frac{f_2^{58}f_3^{21}f_4^5}{f_1^{71}f_6^{12}f_8^2}
- 81 q^2\frac{f_6^5 f_{24}^2}{f_2 f_{12}}   \nonumber \\
& \equiv
- 32 \frac{f_2^{72}f_8^2}{f_4 f_6^{24}}
   \Big(\frac{f_3}{f_1^3}\Big)^{27}(f_1^4)
- 45 \frac{f_2^9 f_8^2}{f_4 f_6^3} \Big(\frac{f_1^3}{f_3}\Big)^3 (f_1^4)
+ 27 \frac{f_4^5}{f_2^2 f_8^2}\Big(\frac{f_3^3}{f_1}\Big)^2
+ 54 \frac{f_2^6 f_{12}^5}{f_6^2 f_{24}^2}f_1^4
     \nonumber \\
& ~~~
- 81 q \frac{f_6^5 f_2 f_8^2}{f_4}\frac{1}{f_1^4}
- 45 q \frac{f_2^{10}f_{24}^2}{f_6^4f_{12}}\Big(\frac{f_1^3}{f_3}\Big)(f_1^4)^2
+ 54 q \frac{f_2^4 f_4^5 f_6^6}{f_8^2}\frac{f_1}{f_3^3}
- 81 q^2 \frac{f_6^5 f_{24}^2}{f_2 f_{12}}
  \pmod{243}.
\end{align}
Using \eqref{f1-4}, \eqref{1-f1-4} and \eqref{f3-f1-3}--\eqref{f1-f3-3} in \eqref{gf-384n+256}, we apply  $E_0$ and Lemma \ref{f-3-r} to deduce 
\begin{align}\label{gf-768n+256}
& \sum_{n=0}^{\infty}\overline{p}(16^2(3n+1))q^n \nonumber \\
& \equiv
- 32 \frac{f_2^{171}f_3^{57}}{f_1^{173} f_4^2 f_6^{54}}
- 45 \frac{f_1^7 f_2^{18}}{f_3^3 f_4^2 f_6^3}
+ 27 \frac{f_2^{11}f_3^4}{f_1^6 f_4^2 f_6^2}
+ 54 \frac{f_1^4 f_2^{10} f_6^5}{f_3^2 f_4^4 f_{12}^2}
- 81 q \frac{f_3^5 f_{12}^2}{f_1 f_6}
- 108 q \frac{f_1^8 f_2^{19} f_6^2 f_{12}^2}{f_3^6 f_4^8} \nonumber \\
& ~~~
- 81 q \frac{f_2 f_3^5 f_4^6}{f_1^9}
+ 117 q \frac{f_1^{10}f_2^{11}f_{12}^2}{f_3^4 f_6^2}
+ \left\{ 81 q \frac{f_1^{13}f_2^2f_4^6f_6}{f_3^5}
- 81 q \frac{f_2^{155}f_3^{55}f_4^6 }{f_1^{167}f_6^{50}} \right\} \nonumber \\
& ~~~
+ \left\{ 27 q \frac{f_2^3 f_6^6}{f_1^2 f_4^2}
- 54 q \frac{f_1^7 f_2^3 f_6^6}{f_3^3 f_4^2} \right\}
- 27 q^2 \frac{f_1^{16}f_4^8f_6^2f_{12}^2}{f_2^5f_3^6}
\nonumber \\
& \equiv
- 32 \frac{f_6^{27}}{f_2^{72}f_4^2}\Big(\frac{f_1^{3}}{f_3}\Big)^{24} \Big(\frac{1}{f_1^2}\Big)
- 45 \frac{f_2^{18}}{f_4^2 f_6^3}\Big(\frac{f_1^3}{f_3}\Big)^3 \Big(\frac{1}{f_1^2}\Big)
+ 27 \frac{f_2^2 f_6}{f_4^2}\Big(\frac{f_1^{3}}{f_3}\Big)(f_3^2) \nonumber \\
& ~~~
+ 54 \frac{f_2^{10}f_6^5}{f_4^4f_{12}^2 }\Big(\frac{f_1^3}{f_3}\Big)^2\Big(\frac{1}{f_1^2}\Big)
- 81 q \frac{f_{12}^2}{f_6}\frac{f_3^3}{f_1}f_3^2
- 108 q \frac{f_2 f_4 f_6^8}{f_{12}} \frac{f_1}{f_3^3}\frac{1}{f_1^2}
- 81 q f_2 f_4^6 (f_3^2)
 \nonumber \\
& ~~~
+ 117 q \frac{f_2^{11} f_{12}^2}{f_6^2}\Big(\frac{f_1^3}{f_3}\Big)^4\Big(\frac{1}{f_1^2}\Big)
- 27 q \frac{f_2^3 f_6^6}{f_4^2}\frac{1}{f_1^2}
- 27 q^2 \frac{f_6^2 f_{12}^5}{f_2^5 f_4} \frac{1}{f_1^2}
  \pmod{243}.
\end{align}
Using \eqref{1-f1-2}, \eqref{f3-2}, and \eqref{f1-3-f3}--\eqref{f1-f3-3}
in \eqref{gf-768n+256}, we apply $E_1$ and Lemma~\ref{f-3-r} to obtain
\begin{align*}
& \sum_{n=0}^{\infty}\overline{p}(16^2(6n+4))q^n \nonumber \\
& \equiv
- 64 \frac{f_2^{72}f_3^{27}f_8^2}{f_1^{77}f_4f_6^{24}}
- 81\left\{
   \frac{f_2^3 f_3^2 f_{12}^5}{f_1^2 f_6 f_{24}^2}
+  \frac{f_1 f_2^6 f_3 f_{12}^5}{f_6^2 f_{24}^2}
+  \frac{f_1^4 f_6 f_{12}^5}{f_2^3 f_{24}^2}  \right\}
-  \left\{
  27 \frac{f_3^6 f_4^5}{f_1^2 f_2^2 f_8^2}
+ 81 \frac{f_1^7 f_3^3 f_4^5}{f_2^2 f_8^2} \right\}
\nonumber \\
&~~~ - 90 \frac{f_1^{13}f_2^9 f_8^2}{f_3^3 f_4 f_6^3}
+ \left\{
  117 \frac{f_2^{66}f_3^{25}f_4^5}{f_1^{75}f_6^{20}f_8^2}
+ 117 \frac{f_1^6 f_2^{12}f_4^5}{f_3^2 f_6^2 f_8^2}
- 108 \frac{f_2^3 f_3 f_4^5 f_6}{f_1^3 f_8^2}
- 81 \frac{f_1^{15}f_2^3 f_4^5 f_6}{f_3^5 f_8^2} \right\} \nonumber \\
&~~~ + 108 \frac{f_1^5 f_2^4 f_3^5 f_8^2}{f_4f_6^4}
+ \left\{
  108 q \frac{f_1^8 f_2^{10} f_6^2 f_8^2}{f_3^4 f_4}
- 54 q \frac{f_2^{64}f_3^{23}f_8^2}{f_1^{73}f_4f_6^{16}}
- 54 q \frac{f_2 f_3^2 f_6^5 f_8^2}{f_1^{10}f_4} \right\}\nonumber \\
&~~~
- 27 q \frac{f_2 f_6^5 f_8^2}{f_1 f_3 f_4}
+ 108 q \frac{f_2^{58}f_3^{21}f_4^5}{f_1^{71}f_6^{12}f_8^2}
- 54 q \frac{f_1^2 f_2 f_3^2 f_{24}^2}{f_6 f_{12}}
- 81 q^2 \frac{f_6^5 f_{24}^2}{f_2 f_{12}}
\nonumber \\
& \equiv
- 64 \frac{f_2^{72}f_3^{27}f_8^2}{f_1^{77}f_4f_6^{24}}
- 108 \frac{f_3^6 f_4^5}{f_1^2 f_2^2 f_8^2}
- 90 \frac{f_1^{13}f_2^9 f_8^2}{f_3^3 f_4 f_6^3}
+ 45 \frac{f_1^6 f_2^{12}f_4^5}{f_3^2 f_6^2 f_8^2}
+ 108 \frac{f_2^4 f_3^8 f_8^2}{f_1^4 f_4 f_6^4}
\nonumber \\
& ~~~
- 27 q \frac{f_2 f_6^5 f_8^2}{f_1 f_3 f_4}
+ 108 q \frac{f_1 f_2^4 f_4^5 f_6^6}{f_3^3 f_8^2}
- 54 q \frac{f_1^2 f_2 f_3^2 f_{24}^2}{f_6 f_{12}}
- 81 q^2 \frac{f_6^5 f_{24}^2}{f_2f_{12}}
  \pmod{243}.
\end{align*}
Hence we get \eqref{gf-16^2(6n+4)}.  
\qed

\begin{proof}[Proof of Theorem~\ref{thm-con-rel}] 
Combining \eqref{gf-6n+4}, \eqref{gf-16(6n+4)}, and \eqref{gf-16^2(6n+4)}, and then applying Lemma \ref{f-3-r},  we find that
\begin{align}\label{gf-con-rel-0}
&\sum_{n=0}^{\infty} \Big(4\cdot\overline{p}(16^2(6n+4))
+ 4\cdot\overline{p}(16(6n+4))
+ \overline{p}(6n+4)\Big) q^n  \nonumber \\
& \equiv
  36 \frac{f_1^{13}f_2^9 f_8^2}{f_3^3 f_4 f_6^3}
- 4 \frac{f_2^{18}f_3^9f_8^2}{f_1^{23}f_4 f_6^6}
- 77 \frac{f_2^{72}f_3^{27} f_8^2}{f_1^{77}f_4f_6^{24}}
- 54 \frac{f_2^4 f_3^8 f_8^2}{f_1^4 f_4 f_6^4}
+ \left\{ 27 q \frac{f_2 f_6^5 f_8^2}{f_1 f_3 f_4}
- 81 q\frac{f_2 f_6^5 f_8^2}{f_1^4f_4} \right\}  \nonumber \\
& ~~~
+ 108 q \frac{f_1 f_2^4 f_4^5 f_6^6}{f_3^3 f_8^2}
- 54 \frac{f_3^6 f_4^5}{f_1^2 f_2^2 f_8^2}
- 90 \frac{f_1^6 f_2^{12}f_4^5}{f_3^2f_6^2  f_8^2}
+ \left\{ 27 q \frac{f_1^2 f_2 f_3^2 f_{24}^2}{f_6 f_{12}}
+ 117 q \frac{f_1^{11}f_2^{10}f_{24}^2}{f_3 f_6^4 f_{12}}  \right\}
\nonumber \\
& ~~~
- \left\{ 81 q^2 \frac{f_2^2 f_6^4 f_{24}^2}{f_{12}}
+ 81 q^2 \frac{f_6^5 f_{24}^2}{f_2 f_{12}}  \right\} \nonumber \\
& \equiv
\left\{  36 \frac{f_1^{13}f_2^9 f_8^2}{f_3^3 f_4 f_6^3}
- 4 \frac{f_2^{18}f_3^9f_8^2}{f_1^{23}f_4 f_6^6}
- 77 \frac{f_2^{72}f_3^{27} f_8^2}{f_1^{77}f_4f_6^{24}}
- 54 \frac{f_2^4 f_3^8 f_8^2}{f_1^4 f_4 f_6^4}
- 54 q \frac{f_2 f_6^5 f_8^2}{f_1 f_3 f_4} \right\}
  \nonumber \\
& \quad
+ \left\{
  108 q \frac{f_1 f_2^4 f_4^5 f_6^6}{f_3^3 f_8^2}
- 54 \frac{f_3^6 f_4^5}{f_1^2 f_2^2 f_8^2}
- 90 \frac{f_1^6 f_2^{12}f_4^5}{f_3^2f_6^2  f_8^2}
\right\}
- 99 q \frac{f_1^{11}f_2^{10}f_{24}^2}{f_3 f_6^4 f_{12}}
+ 81 q^2 \frac{f_6^5 f_{8}^6}{f_2f_{4}^3}  \pmod{243}.
\end{align}
Equations \eqref{f6-3-f2-9} and \eqref{q-f1-f6-6}, together with Lemma \ref{f-3-r}, yield
\begin{align}\label{gf-con-rel-ter-1}
& 36 \frac{f_1^{13}f_2^9 f_8^2}{f_3^3 f_4 f_6^3}
- 4 \frac{f_2^{18}f_3^9f_8^2}{f_1^{23}f_4 f_6^6}
- 77 \frac{f_2^{72}f_3^{27} f_8^2}{f_1^{77}f_4f_6^{24}}
- 54 \frac{f_2^4 f_3^8 f_8^2}{f_1^4 f_4 f_6^4}
 - 54 q \frac{f_2 f_6^5 f_8^2}{f_1 f_3 f_4}   \nonumber \\
& \equiv
36 \frac{f_1^{13}f_2^9 f_8^2}{f_3^3 f_4 f_6^3}
- 4 \frac{f_2^{72}f_3^9f_8^2}{f_1^{23}f_4 f_6^{24}}\left(\frac{f_6^3}{f_2^9}\right)^6
- 77 \frac{f_2^{72}f_3^{27} f_8^2}{f_1^{77}f_4f_6^{24}}
- 54 \frac{f_2^4 f_3^8 f_8^2}{f_1^4 f_4 f_6^4}
- 54 \frac{f_2^3f_3^2f_8^2}{f_1^2 f_4 f_6}
    \left(q \frac{f_1 f_6^6}{f_2^2 f_3^3}\right)  \nonumber \\
& \equiv
27 \frac{f_2^{67}f_3^{32}f_8^2}{f_1^{76}f_4f_6^{25}}
+ \left\{
 36 \frac{f_1^{13}f_2^9 f_8^2}{f_3^3 f_4 f_6^3}
- 54 \frac{f_2^9 f_3^3 f_8^2}{f_1^5 f_4 f_6^3}
+ 135 \frac{f_3^{27}f_2^{72}f_8^2}{f_1^{77}f_6^{24}f_4} \right\}
\nonumber \\
& \equiv  27 \frac{f_3^{8}f_8^2}{f_1^{4}f_2^5f_4f_6}
+ 117 \frac{f_1^{13}f_2^9 f_8^2}{f_3^3 f_4 f_6^3}
 \pmod{243}.
\end{align}
Using \eqref{q-f1-f6-6}, \eqref{psi-phi-Ram-f}, and Lemma \ref{f-3-r}, we obtain 
\begin{align}\label{gf-con-rel-ter-2}
&  108 q \frac{f_1 f_2^4 f_4^5 f_6^6}{f_3^3 f_8^2}
- 54 \frac{f_3^6 f_4^5}{f_1^2 f_2^2 f_8^2}
- 90 \frac{f_1^6 f_2^{12}f_4^5}{f_3^2f_6^2  f_8^2}
\nonumber \\
& \equiv
  108 \frac{f_2^6 f_4^5}{f_8^2}
    \left(q \frac{f_1 f_6^6}{f_2^2 f_3^3}\right)
- 54 \frac{f_3^6 f_4^5}{f_1^2 f_2^2 f_8^2}
- 90 \frac{f_1^6 f_2^{12}f_4^5}{f_3^2f_6^2 f_8^2}
\nonumber \\
& \equiv
- \left\{ 54 \frac{f_3^6 f_4^5}{f_1^2 f_2^2 f_8^2}
+ 108 \frac{f_2^7 f_4^5 f_3^6}{f_8^2 f_1^2 f_6^3} \right\}
+ \left\{
108 \frac{f_2^{12}f_4^5 f_3}{f_8^2 f_1^3 f_6^2}
- 90 \frac{f_1^6 f_2^{12}f_4^5}{f_3^2f_6^2f_8^2} \right\}
\nonumber \\
& \equiv
  81 \frac{f_3^6 f_4^5}{f_1^2 f_2^2 f_8^2}
+ 18 \frac{f_2^{15} f_4^5}{f_6^3 f_8^2}
     \left(\frac{f_1^6 f_6}{f_2^3 f_3^2} \right) \nonumber \\
& \equiv \left\{
  81 \frac{f_3^6 f_4^5}{f_1^2 f_2^2 f_8^2}
+ 162 \frac{f_2^{16} f_3^6 f_4^5}{f_1^2 f_6^6 f_8^2}  \right\}
+ 99 \frac{f_2^{21} f_3 f_4^5}{f_1^3 f_6^5 f_8^2}  \nonumber \\
& \equiv 99 \frac{f_3 f_4^5 f_6^4}{f_1^3 f_2^6 f_8^2}
 \pmod{243}.
\end{align}
Substituting \eqref{gf-con-rel-ter-1} and \eqref{gf-con-rel-ter-2} into \eqref{gf-con-rel-0} gives 
\begin{align}\label{gf-con-rel-1}
&\sum_{n=0}^{\infty} \Big(4\cdot\overline{p}(16^2(6n+4)) + 4\cdot\overline{p}(16(6n+4)) + \overline{p}(6n+4)\Big) q^n  \nonumber \\
& \equiv
27 \frac{f_3^8 f_8^2}{f_1^4 f_2^5 f_4 f_6}
+ 117 \frac{f_1^{13}f_2^9 f_8^2}{f_3^3 f_4 f_6^3}
+ 99 \frac{f_3 f_4^5 f_6^4}{f_1^3 f_2^6 f_8^2}
- 99 q \frac{f_1^{11}f_2^{10}f_{24}^2}{f_3 f_6^4 f_{12}}
+ 81 q^2 \frac{f_6^5 f_{8}^6}{f_2f_{4}^3}
 \pmod{243}.
\end{align}
In view of \eqref{lem-f1-3-f6-2} and Lemma \ref{f-3-r}, it follows that
\begin{align}\label{gf-con-rel-1-1}
117 \frac{f_1^{13}f_2^9 f_8^2}{f_3^3 f_4 f_6^3}
& = 117 \frac{f_1^{19} f_2^6 f_8^2}{f_3^5 f_4 f_6^2}
\left(\frac{f_2^3 f_3^2}{f_1^6 f_6}\right)
\equiv 117 \frac{f_1^{19} f_2^6 f_8^2}{f_3^5 f_4 f_6^2}
\left(\frac{f_1^3 f_6^2}{f_2^6 f_3} + 9 q \frac{f_2^{16}}{f_1^8}\right)  \nonumber \\
& \equiv 117 \frac{f_1^{22}f_8^2}{f_3^6 f_4}
+ 81 q \frac{f_1^{11}f_2^{22}f_8^2}{f_3^5 f_4 f_6^2} \nonumber \\
& \equiv
117 \frac{f_3^3 f_8^2}{f_1^5 f_4}
+ 81 q \frac{f_2^{16} f_8^2}{f_1^4 f_4}
 \pmod{243}.
\end{align}
Utilizing \eqref{lem-f12-2-f4-6}, \eqref{lem-qf4^8-f2^4-m9}, \eqref{f1-3-f3} and Lemma \ref{f-3-r}, we get
\begin{align}\label{gf-con-rel-1-2}
- 99 q \frac{f_1^{11} f_2^{10} f_{24}^2}{f_3 f_6^4 f_{12}}
& \equiv
- 99 q \frac{f_1^{11} f_2^{10} f_8^6}{f_3 f_6^4 f_{12}} \Big(\frac{f_{24}^2}{f_8^6}\Big)  \nonumber \\
& \equiv
- 99 q \frac{f_1^{11} f_2^4 f_8^6}{f_3 f_6^2 f_{12}}
+  54 q \frac{f_1^{11} f_2 f_8^6}{f_3 f_4^3 f_6}  \Big(2 q^2\frac{f_8^8}{f_4^4}\Big)  \nonumber \\
& \equiv
54 q \frac{f_1^{11}f_4^{17}}{f_2^7 f_3 f_6 f_8^2}
+ 90 q \frac{f_1^{11} f_2^4 f_8^6}{f_3 f_4^3 f_6^2} \Big(\frac{f_4^3}{f_{12}}\Big)  \nonumber \\
& \equiv
54 q \frac{f_1^{11} f_4^{17}}{f_2^7 f_3 f_6 f_8^2}
+ 90 q \frac{f_1^{14} f_2^4 f_8^6}{f_3^2 f_4^3 f_6^2}
+ 27 q^2 \frac{f_1^{11} f_2^6 f_{4}^5 f_8^6}{f_3 f_6^4}
 \pmod{243}.
\end{align}
Substituting \eqref{gf-con-rel-1-1} and \eqref{gf-con-rel-1-2} into \eqref{gf-con-rel-1} yields 
\begin{align}\label{gf-con-rel-2}
&\sum_{n=0}^{\infty} \Big(4\cdot\overline{p}(16^2(6n+4)) + 4\cdot\overline{p}(16(6n+4)) + \overline{p}(6n+4)\Big) q^n  \nonumber \\
& \equiv
f_8^2\left\{ 27 \frac{f_3^8}{f_1^4 f_2^5 f_4 f_6}
+ 117 \frac{f_3^3}{f_1^5 f_4}
+ 81 q \frac{f_2^{16} }{f_1^4 f_4} \right\}
\nonumber \\
&
+ f_8^2\left\{ 99 \frac{f_3 f_4^5 f_6^4}{f_1^3 f_2^6 f_8^4}
+ 54 q \frac{f_1^{11} f_4^{17}}{f_2^7 f_3 f_6 f_8^4}
+  90 q \frac{f_1^{14} f_2^4 f_8^4}{f_3^2 f_4^3 f_6^2}
+ 27 q^2 \frac{f_1^{11} f_2^6 f_{4}^5 f_8^4}{f_3 f_6^4}
+ 81 q^2 \frac{f_6^5 f_{8}^4}{f_2f_{4}^3} \right\}   \pmod{243}.
\end{align}
By \eqref{f_4^{10}-f_2^4 f_8^4}, \eqref{8q-f_8^4-f_4^2}, \eqref{gf-con-rel-4-t2}, \eqref{lem-f1-3-f6-2} and \eqref{gf-con-rel-4-t3}, we find that
\begin{align}\label{gf-con-rel-2-ter}
& 99 \frac{f_3 f_4^5 f_6^4}{f_1^3 f_2^6 f_8^4}
+ 54 q \frac{f_1^{11} f_4^{17}}{f_2^7 f_3 f_6 f_8^4}
+  90 q \frac{f_1^{14} f_2^4 f_8^4}{f_3^2 f_4^3 f_6^2}
+ 27 q^2 \frac{f_1^{11} f_2^6 f_{4}^5 f_8^4}{f_3 f_6^4}
+ 81 q^2 \frac{f_6^5 f_{8}^4}{f_2f_{4}^3} \nonumber \\
& \equiv
\left(27q \frac{f_1^{11} f_4^7}{f_2^3f_3f_6}
- 72 \frac{f_3 f_6^4}{f_1^3 f_2^2 f_4^5}
\right)
\left(2 \frac{f_4^{10}}{f_2^4 f_8^4} \right)
+ \left( 72 \frac{f_1^{14}f_2^4}{f_3^2 f_4 f_6^2}
- 27 q \frac{f_1^{11} f_2^6 f_4^7}{f_3 f_6^4}
- 81 q \frac{f_6^5}{f_2 f_4} \right) \left(8 q \frac{f_8^4}{f_4^2}\right)
\nonumber \\
& \equiv  \left\{
72 \frac{f_1^{10}f_2^{14}}{f_3^2 f_4^5 f_6^2}
- 72 \frac{f_1 f_3 f_6^4}{f_2^4 f_4^5}
- 81 q \frac{f_2^{24}}{f_1^4 f_4^5} \right\}
- 72 \frac{f_1^{18}f_2^2}{f_3^2 f_4 f_6^2}
- 72 \frac{f_2^8 f_3 f_6^4}{f_1^7 f_4^9}
+ 54 q \frac{f_1^{15} f_2^4 f_4^7}{f_3f_6^4}
+ 81 q \frac{f_1^4f_2^{12}}{f_4}
 \nonumber \\
& \equiv
- 72 \frac{f_1^{24}}{f_2f_3^4f_4f_6} \left(\frac{f_2^3f_3^2}{f_1^6f_6} \right)
+\left\{ - 72 \frac{f_2^8 f_3 f_6^4}{f_1^7 f_4^9}
+ 54 q \frac{f_1^{15} f_4^7}{f_2^5f_3f_6} \right\}
+ 81 q \frac{f_1^4f_2^{12}}{f_4}
 \nonumber \\
& \equiv
- 36 \frac{f_3^3}{f_1^5 f_4}
- 108 \frac{f_3^{8}}{f_1^4f_2^5 f_4 f_6}
+ 81 q \frac{f_2^{16}}{f_1^4 f_4}
+ 36 \frac{f_3^{4}f_6}{f_2^7 f_4}
- 36 \frac{f_1^9 f_3 f_6^4}{f_2^{16}f_4}
- 81 q \frac{f_1^4f_2^{12}}{f_4}
 \pmod{243}.
\end{align}
Substituting \eqref{gf-con-rel-2-ter} into \eqref{gf-con-rel-2}
and using Lemma \ref{f-3-r} gives
\begin{align*}
&\sum_{n=0}^{\infty} \Big(4\cdot\overline{p}(16^2(6n+4)) + 4\cdot\overline{p}(16(6n+4)) + \overline{p}(6n+4)\Big) q^n  \nonumber \\
& \equiv
f_8^2\left\{
  81 \frac{f_3^3}{f_1^5 f_4}
- 81 \frac{f_3^8}{f_1^4 f_2^5 f_4 f_6}
+ 162 q \frac{f_2^{16} }{f_1^4 f_4}
+ 36 \frac{f_3^{4}f_6}{f_2^7 f_4}
- 36 \frac{f_1^9 f_3 f_6^4}{f_2^{16}f_4}
- 81 q \frac{f_1^{4}f_2^{12}}{f_4} \right\} \nonumber \\
& \equiv
f_8^2\left\{
  81 \frac{f_1^4}{f_4}
- 81 \frac{f_1^{20}}{f_2^8 f_4}
- 81 q \frac{f_2^{16} }{f_1^4 f_4}
+ 36 \frac{f_3^{4}f_6}{f_2^7 f_4}
- 36 \frac{f_1^9 f_3 f_6^4}{f_2^{16}f_4}
- 81 q \frac{f_1^{4}f_2^{12}}{f_4} \right\} \nonumber \\
& \equiv
\frac{f_1^4f_8^2}{f_4}
\left\{
  81 \left(
1 - \frac{f_1^{16}}{f_2^8}
- q \frac{f_2^{16} }{f_1^8} \right)
+ \left( 36 \frac{f_3^{4}f_6}{f_1^4 f_2^7}
- 36 \frac{f_1^5 f_3 f_6^4}{f_2^{16}}
- 81 q f_2^{12} \right) \right\} \nonumber \\
& \equiv 0 \pmod{243}.
\end{align*}
The first bracket vanishes modulo $3$ by \eqref{q-f2^{16}-f1^8}, so its multiple by $81$ vanishes modulo $243$. The second bracket vanishes modulo $243$ by \eqref{81q-f2^{12}}.
This completes the proof. 
\end{proof}

\section{The congruence modulo 1944}
We first establish the required divisibility by $8$.

\begin{lemma}\label{lem-mod8}
For all integers $\alpha,n\ge0$,
\[
\overline p\bigl(16^\alpha(96n+76)\bigr)\equiv0\pmod8.
\]
\end{lemma}

\begin{proof}
Chen, Hou, Sun, and Zhang~\cite[Theorem~1.3]{Chen-Hou-Sun-Zhang-2016}
proved that, for every integer $m\ge0$,
\begin{equation}\label{rel-mod16}
\overline p(4m)\equiv(-1)^m\overline p(m)\pmod{16}.
\end{equation}

Let $S(q)=\sum_{j\ge1}(-1)^j q^{j^2}$. By \eqref{gf-overp}
and \eqref{varphi-q}, 
\begin{equation}\label{gf-mod8}
G(q)=\frac{1}{\varphi(-q)}
=\frac{1}{1+2S(q)}
\equiv 1-2S(q)+4S(q)^2\pmod{8}.
\end{equation}
Every square is congruent to $0$ or $1$ modulo $4$, so neither
a square nor a sum of two squares can be congruent to $3$ modulo $4$.
Since $24n+19\equiv3\pmod4$, the coefficient of $q^{24n+19}$
in \eqref{gf-mod8} vanishes modulo $8$. Hence
$\overline p(24n+19)\equiv0\pmod8$.

Finally,
\[
16^\alpha(96n+76)=4^{2\alpha+1}(24n+19).
\]
The integer $24n+19$ is odd, whereas $4^j(24n+19)$ is even for
$j\ge1$. Thus, iterating \eqref{rel-mod16} gives
\[
\overline p\bigl(16^\alpha(96n+76)\bigr)
\equiv-\overline p(24n+19)
\equiv0\pmod8.
\]
\end{proof}

\begin{proof}[Proof of Theorem~\ref{thm-inf-cong-96-m1944}]
Applying  \eqref{f2f3^6-f1^2f6^3} to the definition of $H(q)$ in \eqref{eq:H},  we obtain
\begin{align*}
H(q)
=\frac{f_2f_6^3}{f_1}
\left(
\frac{f_2f_3^6}{f_1^2f_6^3}
-\frac{f_3f_6}{f_1^3f_2^3}
+q\frac{f_1f_6^6}{f_2^2f_3^3}
\right) 
\equiv 0 \pmod{9}.  
\end{align*}
The final congruence follows from Lemma~\ref{f-3-r} with $k=2$ and $r=2$. 
Since $243\mid 108\cdot9$, Theorem~\ref{thm-ini-gf-96} implies that for every $n\ge0$, 
\begin{equation}\label{cong-initial}
\overline p(96n+76)\equiv
\overline p\bigl(16(96n+76)\bigr)
\equiv0\pmod{243}.
\end{equation}

Fix $n\ge0$ and set
$a_\alpha=\overline p\bigl(16^\alpha(96n+76)\bigr)$
for $\alpha\ge0$. Because
$16^\alpha(96n+76)\equiv4\pmod6$, Theorem~\ref{thm-con-rel}
gives
\[
a_\alpha+4a_{\alpha+1}+4a_{\alpha+2}
\equiv0\pmod{243}.
\]
The initial values $a_0$ and $a_1$ vanish modulo $243$
by \eqref{cong-initial}. As $4$ is invertible modulo $243$,
induction shows that $a_\alpha\equiv0\pmod{243}$ for every
$\alpha\ge0$. Lemma~\ref{lem-mod8} also gives
$a_\alpha\equiv0\pmod8$. Since $\gcd(8,243)=1$, the Chinese
remainder theorem yields $a_\alpha\equiv0\pmod{1944}$.
\end{proof}

\medskip
{\noindent \bf Acknowledgments}
This work was supported by the National Natural Science Foundation of China (Grant No. 11901403) and
the Fundamental Research Funds for Shijiazhuang Tiedao University (Grant No. ZQK202610).

\medskip
{\noindent \bf Data Availability}
No external datasets were used in this article.

\medskip
{\noindent \bf Conflicts of Interest}
The author declares no competing interests.

\end{document}